\documentclass[twoside,12pt]{article}
\title{The extended zero-divisor graphs of the amalgamated duplication of a ring along an ideal}
\date{}
\author{}
\usepackage[all]{xy}
\usepackage{amssymb,amsmath,latexsym}
\usepackage{theorem}
\usepackage{amsfonts}
\usepackage{ulem}
\usepackage{amscd}
\usepackage{amsxtra}
\usepackage{graphicx}
\usepackage{bm}
\usepackage{graphics} %inclusion de figures
\usepackage{graphicx} %inclusion de figures
\usepackage{pstricks,pst-node} %Graphiques
\usepackage{subfigure}
\usepackage{tikz} %Tikz !

\usepackage{tabularx} %Tableaux
\usepackage{multirow} %Gestion des lignes
\usepackage{multicol} %Gestion des colones
\usepackage{arydshln} %Lignes en pointillés
\usepackage{fancybox} %Boites
\usepackage{multicol} %Colonnes
\usepackage{array} %Tableaux maths
\usepackage{fancybox}

\usepackage{mathrsfs}
\usepackage{array,multirow,makecell}
\setcellgapes{1pt}
\makegapedcells

\newcolumntype{R}[1]{>{\raggedleft\arraybackslash }b{#1}}
\newcolumntype{L}[1]{>{\raggedright\arraybackslash }b{#1}}
\newcolumntype{C}[1]{>{\centering\arraybackslash }b{#1}}

\usetikzlibrary{arrows} 

\newcommand{\field}[1]{\mathbb{#1}}

\newcommand{\R}{\field{R}}

\newcommand{\Z }{\field{Z}}
\newcommand{\N }{\field{N}}

\theoremstyle{}\newtheorem{thm}{\bf Theorem}[section]
\theoremstyle{}\newtheorem{cor}[thm]{\bf Corollary}
\theoremstyle{}\newtheorem{lem}[thm]{\bf Lemma}
\theoremstyle{}
\theoremstyle{}
\theoremstyle{}\newtheorem{pro}[thm]{\bf Proposition}
\theoremstyle{}\newtheorem{exm}[thm]{\bf Example}
\theoremstyle{}
\theoremstyle{}
\theoremstyle{}\newcommand{\cqfd}{\hfill$\square$}

\def\pr{{\parindent0pt {\bf Proof.\ }}}

\def\cqfd
{\hspace{1cm}
	\rule{2mm}{2mm}%
	\medbreak%
	\par%
}

\def\ann{{\rm Ann}}
\def\nil{{\rm Nil}}
\begin{document}
	\thispagestyle{empty}
	%%%%%%%%%%%%%%%%%%%%%%%%%%%%%%%%%%%%%%%%%%%%%%%%%%%%%%%%%
	%%%%%%%%%%%%%%%%%%%%%%%%%%%%%%%%%%%%%%%%%%%%%%%%%%%%%%%%%
	
	%%%%%%%%%%%%%%%%%%%%%%%%%%%%%%%%%%%%%%%%%%%%%%%%%%%%%%%%%
	%%%TITLE%%%%%%%%%%%%%%%%%%%%%%%%%%%%%%%%%%%%%%%%%%%%%%%%%
	\maketitle \vspace*{-1.5cm}
	
	%%%%%%%%%%%%%%%%%%%%%%%%%%%%%%%%%%%%%%%%%%%%%%%%%%%%%%%%%
	%%%%%%%%%%%%%%%%%%%%%%%%%%%%%%%%%%%%%%%%%%%%%%%%%%%%%%%%%
	%%%%%%%%%%%%%%%%%%%%%%%%%%%%%%%%%%%%%%%%%%%%%%%%%%%%%%%%%
	%%%NAMES%%%%%%%%%%%%%%%%%%%%%%%%%%%%%%%%%%%%%%%%%%%%%%%%%
	\begin{center}
		{\large\bf  Brahim El Alaoui$^{1, a}$, and  Raja L'hamri$^{1,b}$}
		
		\bigskip
		
		%%%%%%%%%%%%%%%%%%%%%%%%%%%%%%%%%%%%%%%%%%%%%%%%%%%%%%%%%
		%%%ADDRESSES%%%%%%%%%%%%%%%%%%%%%%%%%%%%%%%%%%%%%%%%%%%%%
		
		$^1$   Faculty of Sciences,  Mohammed V University in Rabat,  Morocco.\\
		\noindent $^a$\,brahim$\_$elalaoui2@um5.ac.ma; brahimelalaoui0019@gmail.com; \\   $^b$\,raja.lhamri@um5r.ac.ma; rajaaalhamri@gmail.com  \\[0.2cm]  
	\end{center}
	\bigskip \bigskip
	%%%%%%%%%%%%%%%%%%%%%%%%%%%%%%%%%%%%%%%%%%%%%%%%%%%%%%%%%
	%%%%%%%%%%%%%%%%%%%%%%%%%%%%%%%%%%%%%%%%%%%%%%%%%%%%%%%%%
	%%%%%%%%%%%%%%%%%%%%%%%%%%%%%%%%%%%%%%%%%%%%%%%%%%%%%%%%%
	%%%ABSTRACT%%%%%%%%%%%%%%%%%%%%%%%%%%%%%%%%%%%%%%%%%%%%%%
	\noindent{\large\bf Abstract.}
Let $R$ be a commutative ring and $I$ be an ideal of $R$. The amalgamated duplication of $R$ along $I$ is the subring $R\Join I:=\{(r,r+i)| r\in R, i\in I\}$ of $R\times R$. This paper investigates the extended zero-divisor graph of the amalgamated duplication of $R$ along $I$. The purpose of this work is to study when $\overline{\Gamma}(R\Join I)$ and $\Gamma(R\Join I)$ coincide, to characterize  when $\overline{\Gamma}(R\Join I)$ is complete, and to compute the diameter and the girth of $\overline{\Gamma}(R\Join I)$.\\
	

\small{\noindent{\bf 2020 Mathematics Subject Classification : 13A70, 05C25.}   
		
		%%%%%%%%%%%%%%%%%%%%%%%%%%%%%%%%%%%%%%%%%%%%%%%%%%%%%%%%%
		%%%%%%%%%%%%%%%%%%%%%%%%%%%%%%%%%%%%%%%%%%%%%%%%%%%%%%%%%
		%%%%%%%%%%%%%%%%%%%%%%%%%%%%%%%%%%%%%%%%%%%%%%%%%%%%%%%%%
		%%% Section 1
		%%%%%%%%%%%%%%%%%%%%%%%%%%%%%%%%%%%%%%%%%%%%%%%%%%%%%%%%%
		%%%%%%%%%%%%%%%%%%%%%%%%%%%%%%%%%%%%%%%%%%%%%%%%%%%%%%%%%
		%%%%%%%%%%%%%%%%%%%%%%%%%%%%%%%%%%%%%%%%%%%%%%%%%%%%%%%%%

\section{Introduction}
 Throughout this paper we will assume that $R$ is a commutative ring with identity $1\neq 0$ and $Z(R)$ is the set of zero-divisors of $R$. The idea of associating a graph with a commutative ring goes back to Beck \cite{B88}, where he was mainly interested in colorings. In his work, the vertices of the graph were all elements of the ring. In \cite{AL99}, Anderson and Livingston introduced the zero-divisor graph, denoted $\Gamma(R)$, of $R$, which is a simple graph with vertices $Z(R)^*:=Z(R)\setminus \{0\}$, the set of non-zero zero-divisors of $R$, and for distinct $x, y\in Z(R)^*$, the vertices $x$ and $y$ are adjacent if and only if $xy=0$. They began to study the relationship between ring-theoretic and graph-theoretic properties. In \cite{BMT16}, Bennis et al. gave an extension of the classical zero-divisor graph by introducing the extended zero-divisor graph, denoted by $\overline{\Gamma}(R)$, of $R$, which is a simple graph with the same vertex set as in the classical zero-divisor graph $\Gamma(R)$, and two different vertices $x$ and $y$ are adjacent if and only if $x^ny^m=0$ with $x^n\neq 0$ and $y^m\neq 0$ for some integer $n,m\in \N^*$.\\

For a better understanding of the concept of zero-divisor graphs, several authors have been interested in studying zero-divisor graphs of certain ring constructions (see for example \cite{ACS05, AS06, BMT17, MY08}). In this paper, we are interested in studying the extended zero-divisor graph of the amalgamated duplication of a ring along an ideal. Recall that the amalgamated duplication of $R$ along $I$ is the subring $R\Join I:=\{(r,r+i)| r\in R, i\in I\}$ of $R\times R$.\\
The paper is organized as follows.
In Section 2, we study when the extended zero-divisor graph and the classical zero-divisor graph of the amalgamated duplication of a ring $R$ along an ideal $I$ coincide. In Section 3, we characterize when the extended zero-divisor graph of $R\Join I$ is complete. We also study the diameter and the girth of $\overline{\Gamma}(R\Join I)$.\\

Let us fix some definitions used throughout this paper. Let $x$ be an element of $R$, the annihilator of $x$ is defined as $\ann(x):=\{y\in R| \ xy=0\}$. For an ideal $I$ of $R$, $\sqrt{I}$ means the radical of $I$.  An element $x$ of $R$ is called nilpotent if $x^n=0$ for some positive integers $n$. The set of all nilpotent elements is denoted by $\nil(R):=\sqrt{0}$. A ring $R$ is said to be reduced if $\nil(R)=\{0\}$. The ring $\Z/n\Z$ of residues modulo an integer $n$ is denoted by $\Z_n$.  
For a subset $X$ of $R$, we denote $X^*=X\setminus \{0\}$.
 The set of zero divisors of $R\Join I$ is defined as $Z(R\Join I):=  T_1\cup T_2\cup T_3\cup T_4$ where   $T_1=\{(0,i)| i\in I\}$,  $T_2= \{(i,0)| i\in I\}$,  $T_3= \{(r,r+i)| r\in Z(R)\setminus \{0\},\ i\in I\}$ and   $T_4= \{(r,r+i)| r\in R\setminus Z(R),\ j(r+i)=0 \text{ for some } j\in I\setminus \{0\}\}$.\\

Now let us recall some basic notions of graph theory used in this paper. For a more general background on graph theory, we refer the reader to \cite{B98,CL05}. \\
Let $G=(V,E)$ be a simple (undirected) graph. The graph $G$ is said to be connected if there is a path between any two distinct vertices. Namely, if for each pair of distinct vertices $u$ and $v$ there exists a finite sequence of distinct vertices $u=v_1, \dots, v_n=v$ such that each pair $\{v_i,v_{i+1}\}$ is an edge. The length of such a path is $n-1$ (i.e. the number of edges between $u$ and $v$). For two different vertices $u$ and $v$ in $G$, the distance between $u$ and $v$, denoted by $d(u,v)$, is the length of the shortest path connecting $u$ and $v$, and if no such path exists, we set $d(u,v)=\infty$. The diameter of $G$ is the supremum of the distance between vertices. Namely, the diameter of $G$ is the quantity $diam(G):=\sup\{d(u,v)| u, v \in V\}$. A cycle is a path of length $n\geq 3$ that starts and ends at the same vertex. Namely, it is a path of the form $v_1-v_2-\dots-v_n-v_1$, where $v_i\neq v_j$ if $i\neq j$ for some $n\geq 3$.  The girth of $G$, denoted by $girth(G)$, is the length of a shortest cycle in $G$ if $G$ contains a cycle. But if $G$ does not contain a cycle, then the girth is infinite (i.e., $girth(G)=\infty$).  The graph $G$ is said to be complete if for every two distinct vertices $u$ and $v$ in $V$, $u-v$ is an edge in $E$. \\

\section{When $\overline{\Gamma}(R\Join I)=\Gamma(R\Join I)$?}
In this section we study when $\overline{\Gamma}(R\Join I)$ and $\Gamma(R\Join I)$ coincide. In the general case, it was proved in \cite{BMT16} that $\overline{\Gamma}(R)=\Gamma(R)$ if and only if every nonzero nilpotent element has index $2$ (if $\nil(R)\neq \{0\}$) and $\ann(x)=\ann(x^2)$ for every $x\in Z(R)\setminus \nil(R)$. We give the following lemma which allows us to prove the first main result of this section.

\begin{lem}\label{lem_coincide}
	Let $R$ be a ring and $I$ be an ideal of $R$ such that $I\subseteq Z(R)$ and $\Gamma(R)=\overline{\Gamma}(R)$. Then, $\forall r\in Z(R)\setminus \nil(R)$, $\forall i\in I\setminus \nil(R)$; if $r+i\in \nil(R)$, then   $\ann(r)\subseteq \ann(i)$.
\end{lem}
\pr
Let $r\in Z(R)\setminus \nil(R)$ and $i\in I\setminus \nil(R)$ such that $r+i \in \nil(R)$. Let $\alpha \in \ann(r)$, then $\alpha r=0$. Or, $(r+i)^2=r^2+2ri+i^2=0$ (since $\Gamma(R)=\overline{\Gamma}(R)$). Thus, $\alpha i^2=0$ and so $\alpha \in \ann(i^2)=\ann(i)$.
\cqfd 
The converse of Lemma \ref{lem_coincide} does not hold in general. So, to prove this we give the following counterexample.
 
\begin{exm}
	Consider the ring $R\Join I=\Z_2[X,Y]/<X^3,X^2Y>\Join <\bar{X},\bar{Y}>$.  Then, $I=Z(R)=<\bar{X}, \bar{Y}>$ and $\Gamma(R)\neq \overline{\Gamma}(R)$ since $\bar{X}$ and $\bar{Y}$ are adjacent in $\overline{\Gamma}(R)$ but not in $\Gamma(R)$. By a simple calculus we can easily prove that the necessary condition of Lemma \ref{lem_coincide} holds.
\end{exm}

\begin{thm}\label{thm_coincidence}
	Let $R$ be a ring and $I$ be an ideal of $R$ such that $I\setminus \nil(R)\neq \emptyset$. Then, $\overline{\Gamma}(R\Join I)=\Gamma(R\Join I)$ if and only if  the following assertions hold:
	\begin{enumerate}
		\item $\overline{\Gamma}(R)=\Gamma(R)$.
        \item  $I\subseteq\ann(\nil(R))$.
       % \item  $\forall r\in Z(R)\setminus \nil(R)$, $\forall i\in I\setminus \nil(R)$; if $r+i\in \nil(R)$, then   $\ann(r)\subseteq \ann(i)$.
	\end{enumerate}
\end{thm}
\pr
$\Rightarrow)$ Suppose that $\overline{\Gamma}(R\Join I)=\Gamma(R\Join I)$. Let $r\in \nil(R)$, then $(r,r)\in \nil(R\Join I)$. Thus, $(r,r)^2=0$,  by \cite[Theorem 2.1]{BMT16}. Then, every nilpotent element in $R$ has index $2$. Now, let $r\in Z(R)\setminus \nil(R)$ and $s\in \ann(r^2)$. Then, $sr^2=0$. Thus,  $(s,s),(r,r)\in Z(R\Join I)$ and $(s,s)(r,r)^2=0$. Then, $(s,s)\in \ann((r,r)^2)=\ann((r,r))$,  by \cite[Theorem 2.1]{BMT16}. Then, $s\in \ann(r)$. Hence, $(1)$ holds.\\
Next, let $i\in I$ and $r\in \nil(R)$. Let $j\in I\setminus \nil(R)$. Then, $(r,r+j)\in Z(R\Join I)\setminus \nil(R\Join I)$ (since $r+j\notin \nil(R)$). Since $\overline{\Gamma}(R\Join I)=\Gamma(R\Join I)$, $\ann((r,r+j)^2)=\ann((r,r+j))$. Thus, $(i,0)\in \ann((r,r+j)^2)=\ann((r,r+j))$ and so $ir=0$. \\
%Next, let $r\in Z(R)\setminus \nil(R)$, $i\in I\setminus \nil(R)$ such that $rs=0$ and $r+i\in \nil(R)$. Then, $r^2s=0$ and $(r+i)^2=0$     and so $(s,s)\in \ann((r,r+i)^2)=\ann(r,r+i)$ (since  $\overline{\Gamma}(R\Join I)=\Gamma(R\Join I)$). Thus, $s(r+i)=0$ and hence $(3)$ holds.\\

$\Leftarrow)$ To prove that $\overline{\Gamma}(R\Join I)=\Gamma(R\Join I)$, it is sufficient to prove that every nonzero nilpotent element in $R\Join I$ has index $2$ and for every $(r,r+i)\in Z(R\Join I)\setminus \nil(R\Join I)$, $\ann(r,r+i)=\ann((r,r+i)^2)$. Then, first, let $(r,r+i)\in \nil(R\Join I)=\nil(R)\Join (\nil(R)\cap I)$, then $(r,r+i)^2=(r^2,(r+i)^2)=0$ (Since $r,i\in \nil(R)$ and  $(1)$ holds). Now, let $(r,r+i)\in Z(R\Join I)\setminus \nil(R\Join I)$, we have the following cases:\\
  
\textit{Case 1.} $(r,r+i)=(0,i)\in T_1$. If $i\in I\setminus Z(R)$, then it is clear that $\ann((0,i))=\ann((0,i)^2)$. If $i\in Z(R)\setminus \nil(R)$, then $\ann((0,i)^2)=\{(s,s+j)\in R\Join I| (s+j)i^2=0 \}= \{(s,s+j)\in R\Join I| (s+j)i=0 \}=\ann((0,i))$ (since $\overline{\Gamma}(R)=\Gamma(R)$).\\

\textit{Case 2.} $(r,r+i)\in T_2$. Similar to \textit{Case 1}.\\

\textit{Case 3.} $(r,r+i)\in T_3\setminus \nil(R\Join I)$. Then, we the have following three sub-cases: \\

\textit{Sub-case 1.}  $r\in Z(R)\setminus \nil(R)$ and  $i\in I\setminus \nil(R)$.  If $r+i\notin \nil(R)$, then $\ann((r,r+i)^2)=\{(s,s+j)\in R\Join I| sr^2=0 \text{ and } (s+j)(r+i)^2=0\}=\{(s,s+j)\in R\Join I| sr=0 \text{ and } (s+j)(r+i)=0\}=\ann((r,r+i))$ (since $\overline{\Gamma}(R)=\Gamma(R)$). If $r+i\in \nil(R)$, then $(s,s+j)\in \ann((r,r+i)^2)$ implies that $sr^2=0$ and so $sr=0$, by $(1)$. Then, by  $(2)$ and Lemma $\ref{lem_coincide}$,  $(r+i)(s+j)=0$ which implies that $(s,s+j)\in \ann((r,r+i))$. \\

\textit{Sub-case 2.} $r\in Z(R)\setminus \nil(R)$ and  $i\in I\cap \nil(R)$. Then, $r+i\notin \nil(R)$ and so $\ann((r,r+i)^2)=\ann((r,r+i))$ (see \textit{Sub-case 1}).\\

\textit{Sub-case 3.} $r\in \nil(R)^*$  and $i\in I\setminus \nil(R)$. Let $(s,s+j)\in \ann((r,r+i)^2)$, then $(s+j)(r+i)^2=0$ and so  $(s+j)(r+i)=0$ (since $\overline{\Gamma}(R)=\Gamma(R)$). If $r+i\in Z(R)\setminus\nil(R)$, then   $r+i+(-i)=r\in \nil(R)$ and so by $(2)$ and  Lemma \ref{lem_coincide}, $rs=(r+i+(-i))(s+j+(-j))=0$. Hence, $(s,s+j)\in \ann((r,r+i))$. If $r+i\notin Z(R)$, then $\ann((r,r+i))=\{(s,0)\in R\Join I|\ sr=0, s\in I\}=\{(s,0)\in R\Join I| s\in I\}=\ann((r,r+i)^2)$, by the statements  $(2)$ and Lemma $\ref{lem_coincide}$.\\

\textit{Case 4.} $(r,r+i)\in T_4$. If $r+i\in Z(R)\setminus \nil(R)$. Then, $\ann((r,r+i)^2)=\{(0,j)\in R\Join I| j(r+i)^2=0\}= \{(0,j)\in R\Join I| j(r+i)=0\}= \ann((r,r+i))$, by $(1)$. If $r+i\in \nil(R)$, then $\ann((r,r+i)^2)=\ann(r^2,0)=\{(0,j)\in R\Join I| j\in I \}= \{(0,j)\in R\Join I| j(r+i)=0\}= \ann((r,r+i))$, by  the statements $(2)$ and Lemma $\ref{lem_coincide}$.
\cqfd
In the following, we give an example that verifies the statements of Theorem \ref{thm_complete}. 
\begin{exm}
	Consider the ring  $R\Join I= \R[X,Y]/<X^2,XY>\Join <\bar{X},\bar{Y}>$. We have $Z(R)=<\bar{X},\bar{Y}>=I$ and $\nil(R)=<\bar{X}>$. It is easy to see that the statement $(2)$ of Theorem \ref{thm_coincidence} holds. Also, the statements $(1)$ of Theorem \ref{thm_coincidence} holds. Indeed; Let $r\in Z(R)\setminus \nil(R)$ and $s\in \ann(r^2)$, then $sr^2=0$. Thus, $s\in <\bar{X}>$ and so $rs=0$. Then, $s\in \ann(r)$. Since every nilpotent element has index $2$, $\overline{\Gamma}(R)=\Gamma(R)$ by \cite[Theorem 2.1]{BMT16}.
\end{exm}

To show that the two statements $(1)$ and $(2)$ of Theorem \ref{thm_coincidence} are independent, we give the following examples.

\begin{exm}\label{exm_coincide1}
 Consider the ring $R=\R[X,Y]/<X^2>$ and the ideal $I=<\bar{X},\bar{Y}>$. Then, $I\setminus \nil(R)\neq \emptyset$, and the statement $(1)$ holds since $Z(R)=\nil(R)=<\bar{X}>$ and every non-zero nilpotent element has index $2$. But, the statement $(2)$ does not hold since $\bar{X}\bar{Y}\neq \bar{0}$. On the other hand, the vertices $(\bar{X},\bar{X}+\bar{Y})$ and $(\bar{Y},\bar{0})$ are not adjacent in $\Gamma(\R[X,Y]/<X^2>\Join <\bar{X},\bar{Y}>)$, but they are adjacent in $\overline{\Gamma}(\R[X, Y]/<X^2>\Join <\bar{X},\bar{Y}>)$ (since $(\bar{X},\bar{X}+\bar{Y})^2(\bar{Y},\bar{0})=(\bar{0},\bar{0})$ and $(\bar{X},\bar{X}+\bar{Y})(\bar{Y},\bar{0})\neq(\bar{0},\bar{0})$). Thus,  $\overline{\Gamma}(R\Join I)\neq\Gamma(R\Join I)$.
\end{exm}

\begin{exm}
	Consider the ring $R=\Z_{12}$ and the ideal $I=<\bar{2}>$. Then, $Z(R)=<\bar{2},\bar{3}>$, $\nil(R)=<\bar{6}>$ and $I=\ann(\nil(R))$. But, $\Gamma(R)\neq \overline{\Gamma}(R)$ since $\bar{2}$ is adjacent to $\bar{3}$ in $\overline{\Gamma}(R)$, but not in $\Gamma(R)$. 
\end{exm}

Note that the ring $R\Join I$ is isomorphic to the idealization $R\ltimes I$ if and only if $I$ is a nilpotent ideal of index $2$ in $R$. Then, if $I^2=0$, $\overline{\Gamma}(R\Join I)=\Gamma(R\Join I)$ if and only if $\overline{\Gamma}(R)=\Gamma(R)$. Indeed, assuming that $\overline{\Gamma}(R)=\Gamma(R)$, then by \cite[Theorem 2. 1]{BMT16}, every non-zero nilpotent in $R$ is of index $2$, and for every $x\in Z(R)\setminus \nil(R)$, $\ann(x)=\ann(x^2)$. Thus, we can easily deduce that the conditions $(1)$ and $(4)$ of \cite[Theorem 2.1]{BMT17} hold. For the condition $(3)$ of \cite[Theorem 2.1]{BMT17}, suppose by contrast that there exists $y\in \ann(a)$ for some $a\in Z(R)\setminus \nil(R)$, but $yi\neq 0$ for some $i\in I$. Then, $a+i\in Z(R)\setminus \nil(R)$, since $yi(a+i)=0$ and $I^2=0$. On the other hand, we have $y(a+i)=yi\neq 0$ and $y(a+i)^2=0$, so $\ann(a+i)\neq \ann((a+i)^2)$, a contradiction.\\
For the more general case. Namely, for the case when $I\subseteq \nil(R)$, we have the following theorem.
\begin{thm}\label{thm_coincide_nil}
	Let $R$ be a ring and $I\subseteq\nil(R)$ be an ideal of $R$. Then, $\overline{\Gamma}(R\Join I)=\Gamma(R\Join I)$ if and only if $\overline{\Gamma}(R)=\Gamma(R)$. 
\end{thm}
\pr $\Rightarrow)$ Similar to the proof of Theorem \ref{thm_coincidence}.\\   
 $\Leftarrow)$ Assume that $\overline{\Gamma}(R)=\Gamma(R)$. Let $(r,r+i)\in \nil(R\Join I)$, then $r,r+i\in \nil(R)$. Then, $(r,r+i)^2=(r^2,(r+i)^2)=(0,0)$. Now, let $(r,r+i)\in Z(R\Join I)\setminus \nil(R\Join I)$. Let $(s,s+j)\in \ann((r,r+i)^2)$, then $r^2s=0$ and $(r+i)^2(s+j)=0$. Thus, $sr=0$ and $(s+j)(r+i)=0$ since $r,r+i\notin \nil(R)$ and $\overline{\Gamma}(R)=\Gamma(R)$. Then, $(s,s+j)\in \ann((r,r+i))$.
\cqfd

In the following we give an illustrated example of Theorem \ref{thm_coincide_nil}.
\begin{exm}\label{Exmp_coincide} 
	Consider the ring $R\Join I=\Z_8\Join \{\bar{0}, \bar{4}\}$. 
	Then, the illustrated zero-divisor  and extended zero-divisor graphs of $R\Join I$ are as follows.
	\begin{center}
		\begin{tabular}{cc}
			\hspace*{-1cm} \includegraphics[width=2.25in,height=2.1in]{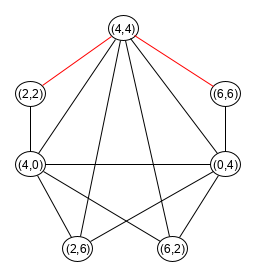} &  \includegraphics[width=2.25in,height=2.1in]{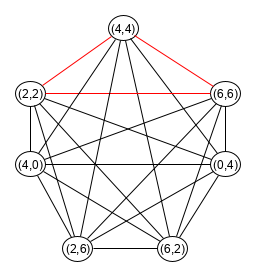}
			\\
			\\
			$\Gamma(\Z_{8}\Join \{\bar{0},\bar{4}\})$  & 
			$\overline{\Gamma}(\Z_{8}\Join \{\bar{0}, \bar{4}\})$
		\end{tabular}
	\end{center}
 We have $I\subseteq \nil(R)$, $\overline{\Gamma}(\Z_{8})\neq \Gamma(\Z_{8})$ (see the red induced subgraphs) and  $\overline{\Gamma}(\Z_{8}\Join \{\bar{0}, \bar{4}\})\neq \Gamma(\Z_{8}\Join \{\bar{0}, \bar{4}\})$.
\end{exm}

%Notice that, if $I=Z(R)$ is an ideal of $R$, we have the following corollary.

%\begin{cor}
%Let $R$ be a ring and $I$ be an ideal of $R$ such that  $I=Z(R)\neq \nil(R)$. Then, $\overline{\Gamma}(R\Join I)=\Gamma(R\Join I)$ if and only if  
%	\begin{itemize}
%		\item $\overline{\Gamma}(R)=\Gamma(R)$, and
%		\item  $\forall r\in Z(R)\setminus \nil(R)$, $\forall i\in I\setminus \nil(R)$; if $r+i\in \nil(R)$ and  $rs=0$ for some $s\in R$, then  $(r+i)(s+j)=0$ for every $j\in I$.
%	\end{itemize} 
%\end{cor}

%\begin{lem}
%	Let $R$ be a ring and $I$ be an ideal of $R$ such that $I^2=0$. Then, $\ann_I(r)=\ann_I(r^2)$ for every $r\in Z(R)\setminus \nil(R)$ if and only if $\ann_I(r+i)=\ann_I((r+i)^2)$ for every $r\in R\setminus Z(R)$ and for every $i\in I$.
%\end{lem}
%\pr $\Rightarrow)$ Let $r\in R\setminus Z(R)$ and  $i\in I$, then $r+i\notin \nil(R)$, otherwise,  $r\in \nil(R)$ (since $i\in \nil(R)$) a contradiction, and so $\ann_I(r+i)=\ann_I((r+i)^2)$.\\
%$\Leftarrow)$ Let  $r\in Z(R)\setminus \nil(R)$. Let  $j\in \ann_I(r^2)$, then $jr^2=jrr=0$. If 
%
%
%$$**************************************$$
%
%
% Next,  let $r\in R\setminus \nil(R)$. Then, $(r,r)\in R\Join I\setminus \nil(R\Join I)$. We may assume that $(r,r)\in Z(R\Join I)\setminus \nil(R\Join I)$. Let $j\in \ann_I(r^2)$, then $(0,j)(r^2,r^2)=0$ and so $(0,j)\in \ann((r,r)^2)=\ann(r,r)$ by \cite[Theorem 2.1]{BMT16}. Thus, $j\in \ann_I(r)$
% 
% 
% $$************************************$$

\section{The diameter and the girth of  $\overline{\Gamma}(R\Join I)$}

In this section, we study the diameter and the girth of the extended zero-divisor graph of the amalgamated duplication $R\Join I$. It was shown in \cite[Theorem 3.1]{BMT16}, that for any commutative ring $R$, the diameter of  $\overline{\Gamma}(R)$ is at most $3$. So, in the following we will study when $diam(\overline{\Gamma}(R\Join I))=1$ (i.e, $\overline{\Gamma}(R\Join I)$ is complete), $diam(\overline{\Gamma}(R\Join I))=2$ and  $diam(\overline{\Gamma}(R\Join I))=3$ (see Theorem \ref{thm_complete}, and Propositions \ref{pro_diameter2} and \ref{pro_diameter3}). We start with the following lemma, which we will need to characterize when $\overline{\Gamma}(R\Join I)$ is complete.

\begin{lem}\label{lem_complete}
	Let $R$ be a ring and $I$ be a non-zero ideal of $R$. Then, the following two statements are equivalent:
	\begin{enumerate}
		\item $Z(R\Join I)=\nil(R\Join I)$ and for every $a,b\in Z(R\Join I)$, $a^{n_a-1}b^{n_b-1}=0$.
		\item $Z(R)=\nil(R)$, $I\subseteq Z(R)$ and for every $x,y\in Z(R)$, $x^{n_x-1}y^{n_y-1}=0$.	
	\end{enumerate}
	 \end{lem}
	 \pr
	 $\Leftarrow)$ $\nil(R\Join I)=\nil(R)\Join (\nil(R)\cap I)= Z(R)\Join (Z(R)\cap I)=Z(R)\Join I$ (since $I\subseteq Z(R)$). Since $I\subseteq \nil(R)$, $T_4=\emptyset$. Thus, $\nil(R\Join I)=Z(R)\Join I=Z(R\Join I)$. Let $a=(r,r+i), b=(s,s+j)\in Z(R\Join I)$, we have $a^{n_a-1}b^{n_b-1}=(r,r+i)^{n_a-1}(s,s+j)^{n_b-1}=(r^{n_a-1}s^{n_b-1}, (r+i)^{n_a-1}(s+j)^{n_b-1})$. Since $n_r\leq n_a$ and $n_s\leq n_b$, $r^{n_a-1}s^{n_b-1}=r^{n_a-n_r}s^{n_b-n_s}r^{n_r-1}s^{n_s-1}=0$ (by (2)). Similarly,  $(r+i)^{n_a-1}(s+j)^{n_b-1}=0$.\\
	 $\Rightarrow)$ Let $r\in Z(R)$, then $(r,r)\in Z(R\Join I)=\nil(R\Join I)$. Then, $r\in \nil(R)$. Thus, $Z(R)=\nil(R)$. Next, let $r\in I$, then $(i,0)\in Z(R\Join I)=\nil(R\Join I)$. Thus, $i\in \nil(R)=Z(R)$. Then, $I\subseteq Z(R)$. Now, let $x,y\in Z(R)$, then $(x,x), (y,y)\in Z(R\Join I)$. Thus, $(x,x)^{n_x-1}(y,y)^{n_y-1}=0$ ($n_{(x,x)}=n_x$ and $n_{(y,y)}=n_y$). Then, $x^{n_x-1}y^{n_y-1}=0$.
	 \cqfd
Note that in the case where $R$ is an integral domain, if $|I|=2$, then the extended zero-divisor graph of $R\Join I$ is just an edge, but if $|I|\geq 3$, then $Z(R\Join I)=T_1\cup T_2$ and so $\overline{\Gamma}(R\Join I)$ is a complete bipartite graph (and then $diam(\overline{\Gamma}(R\Join I))=2$). In the following theorem, we treat the case where $R$ is not an integral domain.
\begin{thm}\label{thm_complete}
 Let $R$ be a ring that is not an integral domain and $I$ be a non-zero ideal of $R$. The following two assertions are equivalent.
	\begin{enumerate}
		\item $\overline{\Gamma}(R\Join I)$ is complete.
		\item $Z(R)=\nil(R)$, $I\subseteq Z(R)$ and for every $x,y\in Z(R)$, $x^{n_x-1}y^{n_y-1}=0$.
		%\item  $Z(R\Join I)=\nil(R\Join I)$ and for every $a,b\in Z(R\Join I)$, $a^{n_a-1}b^{n_b-1}=0$.
	\end{enumerate}
\end{thm}
\pr  $(1)\Rightarrow (2)$. Assume that  $\overline{\Gamma}(R\Join I)$ is complete. Then,  $\overline{\Gamma}(R)$ is complete since  $\overline{\Gamma}(R)$ is an induced subgraph of  $\overline{\Gamma}(R\Join I)$. Thus, by \cite[Theorem 3.3]{BMT16}, $R\cong \Z_2\times \Z_2$ or $Z(R)=\nil(R)$ and for every $x,y\in Z(R)$, $x^{n_x-1}y^{n_y-1}=0$. Suppose that $R\cong \Z_2\times \Z_2$. Then, $Z(R)^*=\{x:=(0,1), y:=(1,0)\}$. We have $1=x+y$. Then, for every $0\neq i\in I$, $i=ix+iy$  and so either $ix\neq 0$ or $iy\neq 0$. Let say $ix\neq 0$, then $(x,x), (0,ix)\in Z(R\Join I)^*$. Thus, $(x,x)$ and $(0,ix)$ are adjacent in $\overline{\Gamma}(R\Join I)$ (since $\overline{\Gamma}(R\Join I)$ is complete). Let $t\in \N^*$ be the smallest positive integer such that $(x,x)^t(0,ix)^s=0$ with $(x,x)^t\neq 0$ and $(0,ix)^s\neq 0$  for some $s\in \N^*$. Then, $i^sx^{t+s}=0$. Since $i=ix+iy$, $i^sx^{s+t-1}=(ix+iy)i^{s-1}x^{s+t-1}=i^sx^{s+t}+i^syx^{s+t-1}=i^sx^{s+t}=0$. Thus, $(x,x)^{t-1}(0,ix)^s=0$ with $(x,x)^{t-1}\neq 0$ and $(0,ix)^s\neq 0$, a contradiction with the fact that $t$ is the smallest positive integer. Hence, $Z(R)=\nil(R)$ and for every $x,y\in R$, $x^{n_x-1} y^{n_y-1}=0$. Next, let $0\neq i\in I$. If $|I|\geq 3$, then for every $i\neq j\in I^*$, $(0,i), (0,j)\in Z(R\Join I)^*$. Then, $(0,i)$ and  $(0,j)$ are adjacent in $\overline{\Gamma}(R\Join I)$. Let $n\in \N^*$ be the smallest positive integer such that $(0,i)^n(0,j)^m=0$ with $(0,i)^n\neq 0$ and $(0,j)^m\neq 0$. Thus, $(j,i)(0,i^{n-1}j^m)=0$ and so $(j,i)\in Z(R\Join I)^*$. Then, $(j,i)^\alpha(0,i)^\beta=0$ with $(j,i)^\alpha \neq 0$ and $(0,i)^\beta \neq 0$ for some $\alpha,\beta \in \N^*$. Thus, $i^{\alpha+\beta}=0$ and so $i\in Z(R)$. Now, if $|I|=2$, then $i^2=0$, otherwise $i^2=i$ and so $i(i-1)=0$. Then,  $i-1\in Z(R)^*$. Thus,  $(1-i,1)\in Z(R\Join I)^*$. Then, for some $\alpha \in \N^*$,  $(0,0)=(1-i,1)^\alpha (0,i)=(0,i)$, a contradiction. Hence, $i\in Z(R)$.\\
$(2)\Rightarrow (1).$ It follows from Lemma \ref{lem_complete}.  
\cqfd 

For example, one can consider $\Z_{8}\Join <\bar{4}>$, $\Z_{p^n}\Join <p>$, and $\R[X]/<X^n>\Join <\bar{X}>$ to get an example of a ring $R\Join I$ that satisfies the assertions of Theorem \ref{thm_complete}.\\

We have the following corollary as a simple consequence of Theorem \ref{thm_complete}.

\begin{cor}\label{cor_complete}
Let $R$ be a ring that is not an integral domain and $I$ be a non-zero ideal of $R$ such that $R\ncong \Z_2\times \Z_2$. Then, $\overline{\Gamma}(R\Join I)$ is complete if and only if $\overline{\Gamma}(R)$ is complete and $I\subseteq Z(R)$.
\end{cor}

 The following example shows that the condition $I\subseteq Z(R)$ in Theorem \ref{thm_complete} and in Corollary \ref{cor_complete} can not be omitted (see also Example \ref{exm_coincide1}).

\begin{exm}
   	      Let $n\geq 2$ be a positive integer and $p$ be a prime number. Consider the ring $R\Join I:=\Z_{p^n}\Join \Z_{p^n}$. Then, $\overline{\Gamma}(\Z_{p^n})$ is complete but $\overline{\Gamma}(\Z_{p^n}\Join \Z_{p^n})$ is not complete since $(1,0)$ is not adjacent to $(p,0)$.
\end{exm}

Now, we give an example that shows that the condition $Z(R)=\nil(R)$ in Theorem \ref{thm_complete} (in particular the condition $\overline{\Gamma}(R)$ is complete in Corollary \ref{cor_complete}) can not be omitted.

\begin{exm}
	Let $n$ be a positive integer and $p$ be a prime number. Consider the ring $R\Join I:= \Z_{p^nq}\Join <p>$.
 We have $<p>\subseteq Z(\Z_{p^nq})$ but $Z(\Z_{p^nq})\neq \nil(\Z_{p^nq})$ (in particular $\overline{\Gamma}(\Z_{p^nq})$ is not complete). On the other hand,   $\overline{\Gamma}(\Z_{p^nq}\Join <p>)$ is not complete (since, for example, $(pq,0)$ and $(q,q)$ are not adjacent in $\overline{\Gamma}(\Z_{p^nq}\Join <p>)$).
\end{exm}

% \begin{exm} 
%    $\overline{\Gamma}(\Z_{p^n}\Join <p>)$ is complete since $<p>\subseteq Z(\Z_{p^n})$ and $\overline{\Gamma}(\Z_{p^n})$ is complete.
%\end{exm}

In \cite[Proposition 4.11]{MY08}, it was shown that $diam(\Gamma(R\Join I))=3$ once $diam(\Gamma(R))=3$. In the following proposition, we give a similar result for the extended zero-divisor graph.

\begin{pro}\label{pro_diameter3}
	Let $R$ be a ring and $I$ be an ideal of $R$. If $diam(\overline{\Gamma}(R))=3$ then, $diam(\overline{\Gamma}(R\Join I))=3$.
\end{pro}
\pr 
Suppose that $diam(\overline{\Gamma}(R\Join I))<3$. Then, $diam(\overline{\Gamma}(R\Join I))=2$. Let $x,y\in Z(R)\setminus\{0\}$ such that $d(x,y)=3$. Thus, $(x,x), (y,y)\in Z(R\Join I)\setminus\{0\}$ and $(x,x)$ and $(y,y)$ are not adjacent in $\overline{\Gamma}(R\Join I)$. Thus, $d((x,x),(y,y))=2$. Then, there exists $(r,r+i)\in Z(R\Join I)\setminus \{0\}$ that is adjacent to both $(x,x)$ and $(y,y)$. Thus, $(r,r+i)^n(x,x)^\alpha=(r,r+i)^n(y,y)^\beta=0$ with $(r,r+i)^n\neq 0$, $(x,x)^\alpha \neq 0$ and $(y,y)^\beta\neq 0$  for some positive integers $n,\alpha,\beta \in \N^*$. If $r^n\neq 0$, then $r$ is adjacent to both $x$ and $y$ in $\overline{\Gamma}(R)$ and so $d(x,y)=2$, a contradiction. Otherwise, $(r+i)^n\neq 0$ and so  $r+i$ is adjacent to both $x$ and $y$ in $\overline{\Gamma}(R)$, then $d(x,y)=2$, a contradiction. Hence, $diam(\overline{\Gamma}(R\Join I))=3$.
\cqfd
The following example shows that the converse of Proposition \ref{pro_diameter3} is not always true.

\begin{exm}
    Consider the ring $R:=\Z_2[X,Y,Z]/<X^3,XY>$ with the ideal $I:=<\bar{X},\bar{Y}, \bar{Z}>$. Then, $d((\bar{Z}, \bar{Z}+(\bar{Y}+\bar{Z})),(\bar{0},\bar{Z}))=3$. This implies that $diam(\overline{\Gamma}(R\Join I))=3$. But $diam(\overline{\Gamma}(R))=2\neq3$.
\end{exm}

In the following result, we characterize when the diameter of $\overline{\Gamma}(R\Join I)$ is equal to $2$ under some assumptions. Let us, first,  say that a graph $G$ has \textit{Condition A} if every edge of $G$ is a part of a triangle. Namely, for every $u,v\in V(G)$ such that $u-v\in E(G)$, there exists a vertex of $G$ that is adjacent to both $u$ and $v$.  

\begin{pro}\label{pro_diameter2}
    Let $R$ be a ring and $I$ be an ideal of $R$ such that $I\subseteq Z(R)$ and $Z(R)$ is an ideal of $R$. If $diam(\overline{\Gamma}(R))=2$ and $\overline{\Gamma}(R)$ has \textit{Condition A}, then $diam(\overline{\Gamma}(R\Join I))=2$.
\end{pro}
\pr Since $diam(\overline{\Gamma}(R))=2$ and $\overline{\Gamma}(R)$ is an induced subgraph of $\overline{\Gamma}(R\Join I)$, $diam(\overline{\Gamma}(R\Join I))\geq 2$. Then, let $(x,x+i)\neq (y,y+j)\in Z(R\Join I)^*$ be two non-adjacent vertices in $\overline{\Gamma}(R\Join I)$. We have the following cases:\\ 
\textit{Case 1.} $x=y=0$. Then, for every $k\in I^*$, $(k,0)$ is adjacent to both $(0,i)$ and $(0,j)$. Thus, $d((0,i),(0,j))= 2$.\\
\textit{Case 2.} $x\neq 0$ and $y=0$. Then, $x\in Z(R)^*$ (since $Z(R)$ is an ideal and $I\subseteq Z(R)$). Thus, there exists $0\neq z\in \ann(x)$. If there exists $k\in I$ such that $zk\neq 0$, then $(zk,0)$ is adjacent to both $(x,x+i)$ and $(0,j)$. Otherwise (i.e, for every $k\in I$, $zk=0$), $(z,z)$ is adjacent to both $(x,x+i)$ and $(0,j)$. Thus, $d((x,x+i),(0,j))=2$.\\
\textit{Case 3.} $x\neq 0$ and $y\neq 0$. Then, $x,y\in Z(R)^*$ (since $Z(R)$ is an ideal and $I\subseteq Z(R)$). If $x$ and $y$ are not adjacent in $\overline{\Gamma}(R)$, then there exists $z\in Z(R)^*$ such that $z$ is adjacent to both $x$ and $y$ (since $diam(\overline{\Gamma}(R))=2$). Then, $z^\alpha x^n=z^\alpha y^m=0$ with $z^\alpha\neq 0$, $x^n\neq 0$ and $y^m\neq 0$ for some positive integers $\alpha,n,m\in \N^*$. If $z^\alpha k\neq 0$ for some $k\in I$, then $(z^\alpha k,0)$ is adjacent to both $(x,x+i)$ and $(y,y+j)$. Otherwise (i.e,  $z^\alpha k=0$ for every $k\in I$), $(z,z)$ is adjacent to both $(x,x+i)$ and $(y,y+j)$. Thus, $d((x,x+i), (y,y+j))=2$. Now, if $x$ and $y$ are adjacent in $\overline{\Gamma}(R)$, then there exists $z\in Z(R)^*$ that is adjacent to both $x$ and $y$ (since $\overline{\Gamma}(R)$ has \textit{Condition A}). Thus,  $z^\alpha x^n=z^\alpha y^m=0$ with $z^\alpha\neq 0$,  $x^n\neq 0$ and $y^m\neq 0$ for some $\alpha,n,m\in \N^*$. Then, by the same argument as above, we prove that $d((x,x+i),(y,y+i))=2$. Therefore, $diam(\overline{\Gamma}(R\Join I))=2$.
\cqfd
The following example shows that the two conditions $diam(\overline{\Gamma}(R))=2$ and $\overline{\Gamma}(R)$ has Condition A, are independent.

\begin{exm}
    Consider the rings $R:=\Z_p\times\Z_q$ and $S:=\Z_{p^n}$, where $p$ and $q$ are two prime numbers and $n\geq 2$ be a positive integer. Then, $diam(\overline{\Gamma}(R))=2$, but $\overline{\Gamma}(R)$ does not verify Condition A. And, $\overline{\Gamma}(S)$ has Condition A, but $diam(\overline{\Gamma}(S))=1$.    
\end{exm}

As an example of a ring $R$ such that $diam(\overline{\Gamma}(R))=2$ and $\overline{\Gamma}(R)$ has Condition A, we can consider the ring $R=\R[X,Y,Z]/<XY,XZ,YZ,X^3>$.\\

%The diameter of the extended zero-divisor graph of a non-reduced ring $R$ was studied when $Z(R)=\nil(R)$ in \cite[Theorem 3.5]{BMT16}. The next proposition gives a similar result for the amalgamated duplication of $R$ along an ideal $I$. For this, we set 
%$\overline{Z}(R):=\{x^{n_x-1}| x\in \nil(R)^*\}$ and $\overline{Z}(R)^2=\{x^{n_x-1}y^{n_y-1}| x,y\in \nil(R)\}$.

%\begin{pro}\label{pro_diameterleq2}
%	Let $R$ be a  ring that is not an integral domain and  $I\subseteq Z(R)=\nil(R)\neq \{0\}$ be a non-zero ideal of $R$. Then, $diam(\overline{\Gamma}(R\Join I))\leq 2$ and exactly one of the following two cases must occur. 
%	\begin{enumerate}
 %           \item If $\overline{Z}(R)^2=\{0\}$, then $diam(\overline{\Gamma}(R\Join  I))=1$. 
%		\item   If $\overline{Z}(R)^2\neq\{0\}$, then $diam(\overline{\Gamma}(R\Join  I))=2$.
%	\end{enumerate}  
%\end{pro}
%\pr  
 %If $\overline{Z}(R)^2=\{0\}$, then by Theorem \ref{thm_complete}, $\overline{\Gamma}(R\Join I)$ is complete and so $diam(\overline{\Gamma}(R\Join I))=1$. Otherwise, there exist $x,y\in Z(R)^*$ such that $x^{n_x-1}y^{n_y-1}\neq 0$. Then, $xy\notin \{0,x,y\}$ and $x-xy-y$ is a path in $\overline{\Gamma}(R)$. This implies that $diam(\overline{\Gamma}(R))=2$. Since $x^{n_x-1}y^{n_y-1}\in \sqrt{\ann(x^n,y^m)\setminus\{0\}}$ for every $n,m\in \N^*$, $diam(\overline{\Gamma}(R\Join I))=2$, by Proposition \ref{pro_diameter2}.  
%\cqfd 

 Now, we will characterize the girth of $\overline{\Gamma}(R\Join I)$.
 %If $|I|=1$, then $\overline{\Gamma}(R)=\overline{\Gamma}(R\Join I)$ and so $girth(\overline{\Gamma}(R))=girth(\overline{\Gamma}(R\Join I))$. So in the following, we will only consider the case where $|I|\geq2$.
 In the case where $R$ is an integral domain, if $|I|=2$ then $\overline{\Gamma}(R\Join I)$ is an edge, and so $girth(\overline{\Gamma}(R\Join I))=\infty$. If $|I|\geq 3$, then $Z(R\Join I)=T_1\cup T_2$. Thus,   $\overline{\Gamma}(R\Join I)$ is a complete bipartite graph and hence $girth(\overline{\Gamma}(R\Join I))=4$.
 %is reduced and thus $\overline{\Gamma}(R\Join I)=\Gamma(R\Join I)$. So, by \cite[Corollary 3.3]{MY08}, $girth(\overline{\Gamma}(R\Join I))=4$ if and only if $R$ is an integral domain and $|I|\geq 3$, and $girth(\overline{\Gamma}(R\Join I))=\infty$ if and only if $I=R\cong \Z_2$. 
 In the following theorem, we consider only the case where $R$ is not an integral domain.
\begin{thm}
	Let $R$ be a ring and $I$ be a non-zero ideal of $R$. If $R$ is not an integral domain, then  $girth(\overline{\Gamma}(R\Join I))=3$. 
\end{thm}
\pr Suppose that $R$ is not an integral domain. Since the classical zero-divisor graph  $\Gamma(R\Join I)$ is a subgraph of $\overline{\Gamma}(R\Join I)$, $girth(\overline{\Gamma}(R\Join I))\leq girth(\Gamma(R\Join I))$. On the other hand, the girth of $\Gamma(R\Join I)$ is $3$,  by \cite[Proposition 3.1]{MY08}. Thus, $girth(\overline{\Gamma}(R\Join I))=3$.\\ 
\cqfd

\begin{cor}
    Let $R$ be a ring and $I$ be a non-zero ideal of $R$. Then, we have the following equivalents: 
    \begin{enumerate}
        \item $girth(\overline{\Gamma}(R\Join I))=3$ if and only if $R$ is not an integral domain.
        \item  $girth(\overline{\Gamma}(R\Join I))=4$ if and only if $R$ is an integral domain and $|I|\geq 3$.
        \item $girth(\overline{\Gamma}(R\Join I))=\infty$ if and only if $R$ is an integral domain and $|I|=2$.
    \end{enumerate}
\end{cor}

%\section{When $\overline{\Gamma}(R\Join I)$ is complemented?}

\end{document}